\documentclass[11pt,a4paper]{article}
\usepackage{mathrsfs}
\usepackage{epsf,epsfig,amsfonts,amsgen,amsmath,amstext,amsbsy,amsopn,amsthm,lineno}
\usepackage{color}
\newtheorem{theorem}{Theorem}

\newtheorem{prop}{Proposition}
\newtheorem{lemma}{Lemma}
\newtheorem{false statement}{False statement}

\theoremstyle{definition}

\newtheorem{claim}{Claim}

\begin{document}

\title{\bf\Large Covering the edges of digraphs in $\mathscr{D}(3,3)$ and $\mathscr{D}(4,4)$ with directed cuts \thanks{Supported  by NSFC
(No.10871158 and 11171273) and Foundation of Shaanxi Educational Committee (No.09JK609)}}

\date{}

\author{Yandong Bai, Binlong Li, Shenggui Zhang\thanks{Corresponding author. E-mail address: sgzhang@nwpu.edu.cn (S. Zhang).}    \\[2mm]
\small Department of Applied Mathematics,  \small Northwestern Polytechnical University,\\
\small Xi'an, Shaanxi 710072, P.R.~China\\} \maketitle

\begin{abstract}
For nonnegative integers $k$ and $l$, let $\mathscr{D}(k,l)$ denote the family of digraphs
in which every vertex has either indegree at most $k$ or outdegree
at most $l$. In this paper we prove that the edges of every digraph
in $\mathscr{D}(3,3)$ and $\mathscr{D}(4,4)$ can be covered by at most five directed
cuts and present an example in $\mathscr{D}(3,3)$ showing that this
result is best possible.

\smallskip
\end{abstract}

\section{Introduction}

In this paper we only consider directed graphs, called here {\em
digraphs}, without loops and parallel edges. We use Bondy and Murty
\cite{bondy_murty} for {terminology and notation} not defined here.

Let $D$ be a digraph. The vertex set and edge set of $D$ are denoted
by $V(D)$ and $E(D)$, respectively. For a vertex $v$ of $D$,  {its
{\em indegree} and {\em outdegree} are denoted by $d_D^-(v)$ and $d_D^+(v)$,
respectively, and its {\em degree} is $d_D(v)=d_D^-(v)+d_D^+(v)$. }
For a bipartition $(X,Y)$ of $V(D)$,
the edge set $E(X,Y)=\{xy\in E(D)\colon\, x\in X, y\in Y\}$
 {is} the {\em directed cut} induced by $(X,Y)$.

 {Harary {\em et al.} \cite{Harary et al.} considered the problem of
covering the edges of an undirected graph with bipartite subgraphs.
They proved that the minimum number of bipartite subgraphs required to cover the
edges of an undirected graph $G$ is $\lceil\log_{2}\chi(G)\rceil$, where $\chi(G)$ is the
chromatic number of $G$. Alon {\em et al.} \cite{Alon et al.} discussed
the problem of covering the edges of a digraph with} {directed cuts}.
They first considered this problem for complete digraphs,
in which every pair of vertices  induces two edges, one in each direction.
 {For convenience, we abbreviate 'directed cut' to 'cut' in the following.}

\begin{theorem}[Alon \emph{et al.} \cite{Alon et al.}]
The minimum number of  {cuts} required to cover the edges of
the complete digraph on $n$ vertices is equal to  {$c(n)$, where}
$$
c(n) {=}\min \{k\colon\, \binom {k}{\lfloor k/2\rfloor}\geq
n\}=\log_2{n}+\frac {1}{2}\log_2{\log_2{n}}+O(1).
$$
\end{theorem}

It follows immediately from Theorem 1 that every digraph whose
underlying graph has chromatic number at most $n$ can be covered by
$c(n)$  {cuts}, as we can use a coloring of the underlying
graph to group the vertices of our digraph into $n$ classes.
 {For nonnegative integers $k$ and $l$}, let $\mathscr{D}(k,l)$ denote the family of digraphs in
which every vertex has either indegree at most $k$ or outdegree at
most $l$. Alon {\em et al.} \cite{Alon et al.} showed that the
underlying graph of every digraph in $\mathscr{D}(k,l)$  {has chromatic number at most $2k+2l+2$}.
This implies the following result.

\begin{theorem}[Alon \emph{et al.} \cite{Alon et al.}]
Every digraph in $\mathscr{D}(k,l)$ can be covered by at most $c(2k+2l+2)$
 {cuts}.
\end{theorem}

By Theorem 2, every digraph in $\mathscr{D}(k,k)$ can be covered by at most
$c(4k+2)$  {cuts}. Here we give  {a better bound}.

\begin{theorem}
Every digraph in $\mathscr{D}(k,k)$ can be covered by at most $c(2k+1)+1$
 {cuts}.
\end{theorem}

\begin{proof}
Let $D\in \mathscr{D}(k,k)$ {, and let} $(X,Y)$ be a bipartition of $V(D)$ such that
$d_D^-(x)\leq k$ for every $x\in X$ and $d_D^+(y)\leq k$ for every
$y\in Y$. Let $D'$ be an arbitrary subdigraph of  {$D-E(X,Y)$}.
Since $xy\notin E(D')$ for any $x\in X$ and $y\in Y$,
 {all the edges of $D'$ leaving the vertices of $X$ are
counted among the edges of $D'$ entering $X$, and hence}
$$
\sum_{x\in V(D')\cap X}d^+_{D'}(x)\leq \sum_{x\in V(D')\cap
X}d^-_{D'}(x)\leq k|V(D')\cap X|.
$$
Similarly we have
$$
\sum_{y\in V(D')\cap Y}d^-_{D'}(y)\leq \sum_{y\in V(D')\cap
Y}d^+_{D'}(y)\leq k|V(D')\cap Y|.
$$
Therefore,
$$
\sum_{v\in V(D')}d_{D'}(v)\leq 2k|V(D')\cap X|+2k|V(D')\cap
Y|=2k|V(D')|.
$$
This implies that $D'$ contains a vertex with degree at most $2k$.
Thus,  {$D-E(X,Y)$} is $2k$-degenerate, and hence its underlying
graph has chromatic number at most $2k+1$  {(see \cite{Szekeres et al.})}.
Thus {$D-E(X,Y)$} can be covered by $c(2k+1)$  {cuts}.
 {With $E(X,Y)$, these cuts cover $D$.}
\end{proof}

 {The bound in Theorem 3 is} not tight. For $k=1$, Theorem 3
implies that every digraph in $\mathscr{D}(1,1)$ can be covered by at most
four  {cuts}, whereas Alon {\em et al.} \cite{Alon et al.}
proved that three  {cuts} suffice. For $k=2$, Theorem 3
implies that every digraph in $\mathscr{D}(2,2)$ can be covered by at most
five  {cuts}, whereas it was  {noted} in
\cite{Lehel et al.} that Rizzi had proved that four  {cuts}
suffice. Examples from \cite{Alon et al.,Lehel et al.} show
that these bounds are best possible.

In this paper we consider  {an} improvement of Theorem 3 for the cases
$k=3$ and $k=4$. From Theorem 3, we know that every digraph in $\mathscr{D}(3,3)$
and $\mathscr{D}(4,4)$ can be covered by at most six  {cuts}. Here we
prove that five  {cuts} suffice.

\begin{theorem}
Every digraph in $\mathscr{D}(3,3)$ and $\mathscr{D}(4,4)$ can be covered by at most
five  {cuts}.
\end{theorem}

In Section 2 we show that the result of Theorem 4 is best possible
by constructing a digraph in $\mathscr{D}(3,3)$ that cannot be covered by four
 {cuts}. In Section 3 we establish a number of results that
will be used in the proof of Theorem 4.  {We complete the proof of}
Theorem 4 in Section 4.

\section{A digraph in $\mathscr{D}(3,3)$ that cannot be covered by four  {cuts}}

 {Let $D_1$ be the orientation of $K_7$ with vertices
$x_1,\ldots,x_7$ such that the out-neighbors of each vertex are the
vertices with the three next higher indices (modulo 7).
Let $D_2$ be a copy of $D_1$ with vertices $y_1,\ldots,y_7$ in order.
Now give each set $\{x_i,x_j,x_k\}$ of three distinct vertices in $D_1$ a
common out-neighbor $z_{i,j,k}$, and let $Z$ be this set of 35 vertices.
Finally, add an edge from each vertex of $V(D_1)\cup Z$ to each vertex of $D_2$.
Denote the resulting digraph  by $D^*$ {(See Fig. 1)}.
It is easy to see that $d^-_{D^*}(x)=d^+_{D^*}(y)=d^-_{D^*}(z)=3$ for $x\in
V(D_1)$, $y\in V(D_2)$ and $z\in Z$, so $D^*\in \mathscr{D}(3,3)$. }

\begin{center}
\setlength{\unitlength}{0.75pt}
\begin{picture}(420,540)
\thicklines

\put(120,420){\put(0,100){\circle*{6}} \put(78,62){\circle*{6}} \put(97,-22){\circle*{6}} \put(43,-90){\circle*{6}} \put(-43,-90){\circle*{6}} \put(-97,-22){\circle*{6}} \put(-78,62){\circle*{6}}
\put(3,103){$x_1$} \put(81,65){$x_2$} \put(100,-19){$x_3$} \put(45,-97){$x_4$} \put(-55,-100){$x_5$} \put(-109,-32){$x_6$} \put(-90,68){$x_7$}
\qbezier(0,100)(0,100)(78,62) \qbezier(0,100)(0,100)(97,-22) \qbezier(0,100)(0,100)(43,-90)
\qbezier(78,62)(78,62)(97,-22) \qbezier(78,62)(78,62)(43,-90) \qbezier(78,62)(78,62)(-43,-90)
\qbezier(97,-22)(97,-22)(43,-90) \qbezier(97,-22)(97,-22)(-43,-90) \qbezier(97,-22)(97,-22)(-97,-22)
\qbezier(43,-90)(43,-90)(-43,-90) \qbezier(43,-90)(43,-90)(-97,-22) \qbezier(43,-90)(43,-90)(-78,62)
\qbezier(-43,-90)(-43,-90)(-97,-22) \qbezier(-43,-90)(-43,-90)(-78,62) \qbezier(-43,-90)(-43,-90)(0,100)
\qbezier(-97,-22)(-97,-22)(-78,62) \qbezier(-97,-22)(-97,-22)(0,100) \qbezier(-97,-22)(-97,-22)(78,62)
\qbezier(-78,62)(-78,62)(0,100) \qbezier(-78,62)(-78,62)(78,62) \qbezier(-78,62)(-78,62)(97,-22)
\put(0,100){\put(-18,-9){\vector(2,1){0}} \put(-12,-16){\vector(3,4){0}} \put(-4,-19){\vector(1,4){0}}}
\put(78,62){\put(-18,9){\vector(2,-1){0}} \put(-20,0){\vector(1,0){0}} \put(-18,-9){\vector(2,1){0}}}
\put(97,-22){\put(-4,19){\vector(1,-4){0}} \put(-12,16){\vector(3,-4){0}} \put(-18,9){\vector(2,-1){0}}}
\put(43,-90){\put(12,16){\vector(-3,-4){0}} \put(4,19){\vector(-1,-4){0}} \put(-4,19){\vector(1,-4){0}}}
\put(-43,-90){\put(20,0){\vector(-1,0){0}} \put(18,9){\vector(-2,-1){0}} \put(12,16){\vector(-3,-4){0}}}
\put(-97,-22){\put(12,-16){\vector(-3,4){0}} \put(18,-9){\vector(-2,1){0}} \put(20,0){\vector(-1,0){0}}}
\put(-78,62){\put(-4,-19){\vector(1,4){0}} \put(4,-19){\vector(-1,4){0}} \put(12,-16){\vector(-3,4){0}}}
\qbezier[10](120,0)(120,50)(85,85) \qbezier[10](85,85)(50,120)(0,120) \qbezier[10](0,120)(-50,120)(-85,85) \qbezier[10](-85,85)(-120,50)(-120,0)
\qbezier[10](-120,0)(-120,-50)(-85,-85) \qbezier[10](-85,-85)(-50,-120)(0,-120) \qbezier[10](0,-120)(50,-120)(85,-85) \qbezier[10](85,-85)(120,-50)(120,0)}

\put(120,120){\put(0,100){\circle*{6}} \put(78,62){\circle*{6}} \put(97,-22){\circle*{6}} \put(43,-90){\circle*{6}} \put(-43,-90){\circle*{6}} \put(-97,-22){\circle*{6}} \put(-78,62){\circle*{6}}
\put(3,103){$y_1$} \put(81,65){$y_2$} \put(100,-19){$y_3$} \put(45,-97){$y_4$} \put(-55,-100){$y_5$} \put(-109,-32){$y_6$} \put(-90,68){$y_7$}
\qbezier(0,100)(0,100)(78,62) \qbezier(0,100)(0,100)(97,-22) \qbezier(0,100)(0,100)(43,-90)
\qbezier(78,62)(78,62)(97,-22) \qbezier(78,62)(78,62)(43,-90) \qbezier(78,62)(78,62)(-43,-90)
\qbezier(97,-22)(97,-22)(43,-90) \qbezier(97,-22)(97,-22)(-43,-90) \qbezier(97,-22)(97,-22)(-97,-22)
\qbezier(43,-90)(43,-90)(-43,-90) \qbezier(43,-90)(43,-90)(-97,-22) \qbezier(43,-90)(43,-90)(-78,62)
\qbezier(-43,-90)(-43,-90)(-97,-22) \qbezier(-43,-90)(-43,-90)(-78,62) \qbezier(-43,-90)(-43,-90)(0,100)
\qbezier(-97,-22)(-97,-22)(-78,62) \qbezier(-97,-22)(-97,-22)(0,100) \qbezier(-97,-22)(-97,-22)(78,62)
\qbezier(-78,62)(-78,62)(0,100) \qbezier(-78,62)(-78,62)(78,62) \qbezier(-78,62)(-78,62)(97,-22)
\put(0,100){\put(-18,-9){\vector(2,1){0}} \put(-12,-16){\vector(3,4){0}} \put(-4,-19){\vector(1,4){0}}}
\put(78,62){\put(-18,9){\vector(2,-1){0}} \put(-20,0){\vector(1,0){0}} \put(-18,-9){\vector(2,1){0}}}
\put(97,-22){\put(-4,19){\vector(1,-4){0}} \put(-12,16){\vector(3,-4){0}} \put(-18,9){\vector(2,-1){0}}}
\put(43,-90){\put(12,16){\vector(-3,-4){0}} \put(4,19){\vector(-1,-4){0}} \put(-4,19){\vector(1,-4){0}}}
\put(-43,-90){\put(20,0){\vector(-1,0){0}} \put(18,9){\vector(-2,-1){0}} \put(12,16){\vector(-3,-4){0}}}
\put(-97,-22){\put(12,-16){\vector(-3,4){0}} \put(18,-9){\vector(-2,1){0}} \put(20,0){\vector(-1,0){0}}}
\put(-78,62){\put(-4,-19){\vector(1,4){0}} \put(4,-19){\vector(-1,4){0}} \put(12,-16){\vector(-3,4){0}}}
\qbezier[10](120,0)(120,50)(85,85) \qbezier[10](85,85)(50,120)(0,120) \qbezier[10](0,120)(-50,120)(-85,85) \qbezier[10](-85,85)(-120,50)(-120,0)
\qbezier[10](-120,0)(-120,-50)(-85,-85) \qbezier[10](-85,-85)(-50,-120)(0,-120) \qbezier[10](0,-120)(50,-120)(85,-85) \qbezier[10](85,-85)(120,-50)(120,0)}

\put(320,270){\put(0,100){\put(0,0){\circle*{6}} \put(3,-3){$z_{1,2,3}$} \qbezier(0,0)(0,0)(-200,150) \qbezier(0,0)(0,0)(-122,112) \qbezier(0,0)(0,0)(-103,28)
\put(-50,38){\vector(4,-3){0}} \put(-50,46){\vector(1,-1){0}} \put(-50,14){\vector(4,-1){0}}}
\put(0,-50){\put(0,0){\circle*{6}} \put(3,-3){$z_{5,6,7}$} \qbezier(0,0)(0,0)(-243,110) \qbezier(0,0)(0,0)(-297,178) \qbezier(0,0)(0,0)(-278,262)
\put(-50,23){\vector(2,-1){0}} \put(-50,30){\vector(3,-2){0}} \put(-50,47){\vector(1,-1){0}}}
\multiput(0,0)(0,25){3}{\circle*{2}}
\put(0,25){\qbezier[25](20,0)(0,200)(-20,0) \qbezier[25](20,0)(0,-200)(-20,0)
\put(20,-5){$\left.\begin{array}{l} \\ \\ \\ \\ \\ \\ \end{array}\right\}$ 35 vertices}}}

\put(120,270){\put(-60,30){\line(0,-1){60}} \put(60,30){\line(0,-1){60}} \put(-30,20){\line(0,-1){40}} \put(30,20){\line(0,-1){40}}
\multiput(-10,0)(10,0){3}{\circle*{2}}
\put(-60,-3){\vector(0,-1){0}} \put(-30,-3){\vector(0,-1){0}} \put(30,-3){\vector(0,-1){0}} \put(60,-3){\vector(0,-1){0}}}

\put(245,210){\put(-20,60){\line(-1,-1){40}} \put(-10,30){\line(-1,-1){20}} \put(30,-10){\line(-1,-1){20}} \put(60,-20){\line(-1,-1){40}}
\multiput(-8,8)(8,-8){3}{\circle*{2}}
\put(-42,38){\vector(-1,-1){0}} \put(-22,18){\vector(-1,-1){0}} \put(18,-22){\vector(-1,-1){0}} \put(38,-42){\vector(-1,-1){0}}}

\end{picture}

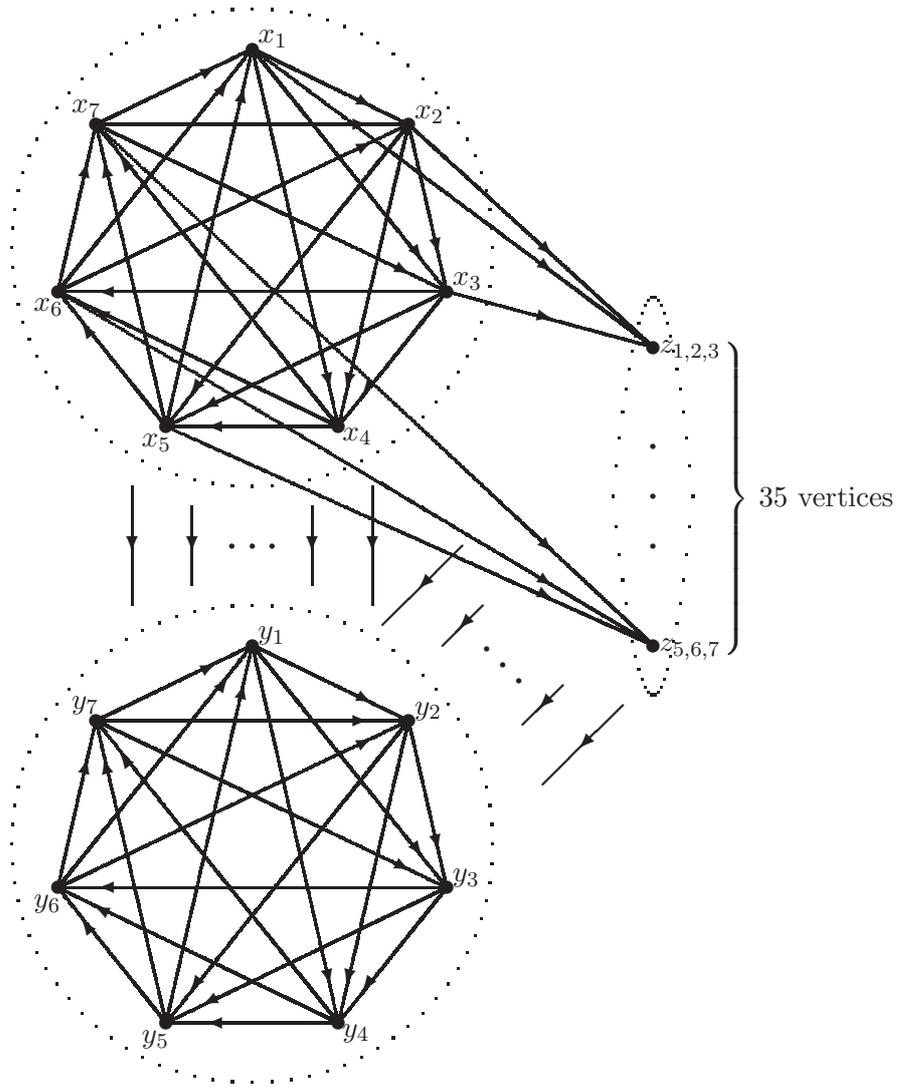
\small Fig. 1. The digraph {$D^*$}.
\end{center}

Now we show that $D^*$ cannot be covered by four  {cuts}.

Assume to the contrary that  {there exist cuts
$E(A_1,B_1),\ldots,E(A_4,B_4)$ that cover $E(D^*)$}. For
$v\in V(D^*)$, let $A(v)=\{i\colon\, v\in A_i\}$ and $B(v)=\{j\colon\, v\in
B_j\}$ {, so} $|A(v)|+|B(v)|=4$.

\begin{claim}
If $uv\in E(D^*)$, then $A(u)\neq A(v)$.
\end{claim}

\begin{proof}
 {It is immediate that $A(u)=A(v)$ prevents $uv$ from being covered.}
\end{proof}

\begin{claim}
 {Neither $D_1$ nor $D_2$ can be covered by three cuts.}
\end{claim}

\begin{proof}
 {
In a regular digraph, any bipartition has the same number
of edges in each direction, and in an orientation of $K_7$ at most 12 edges
cross any bipartition. Hence three cuts cover at most 18 edges,
but $D_1$ and $D_2$ have 21 edges.
}
\end{proof}

\begin{claim}
 {
$|A(x)|\in \{2,3\}$ for $x\in V(D_1)$, and
$|A(y)|\in \{1,2\}$ for $y\in V(D_2)$.
}
\end{claim}

\begin{proof}
It follows from $d^-_{D^*}(x)\geq 1$ and $d^+_{D^*}(x)\geq 1$ that
$1\leq |A(x)|\leq 3$.  {If $A(x)=\{j\}$ for some $x\in V(D_1)$ and index $j$,
then $V(D_2)\subseteq B_j$, since each vertex of $D_2$ is an out-neighbor of $x$.
This requires $D_2$ to be covered by the other three cuts, contradicting Claim 2.
The second assertion can be proved similarly.}
\end{proof}

 {
For $v\in V(D_1)\cup V(D_2)$, it follows that $1\leq |A(v)|\leq 3$.
Since $\{1,2,3,4\}$ has 14 nonempty proper subsets, by Claim 1 each
nonempty subset of $\{1,2,3,4\}$ occurs as $A(v)$ for exactly
one vertex $v\in V(D_1)\cup V(D_2)$. By Claim 3, the vertices of $D_1$
correspond to all four 3-sets and three 2-sets, while those of $D_2$
correspond to all four 1-sets and three 2-sets}.
In the following we use $x_p,x_q,x_r$ to denote the
three vertices in $D_1$ with $|A(x_p)|=|A(x_q)|=|A(x_r)|=2$.

\begin{claim}
$A(x_p)\cap A(x_q)\cap A(x_r)=\emptyset$.
\end{claim}

\begin{proof}
Let $U=\{x\in V(D_1)\colon\, |A(x)|=3\}$ {, so} $|U|=4$. Since $|V(D_1)|=7$,
 {some two consecutive vertices (modulo 7) are in $U$.}
We assume without loss of generality that $x_1,x_2\in U$ and $B(x_1)=\{1\}, B(x_2)=\{2\}$.

If $|B(x_3)|=1$, then we can assume without loss of generality that
$B(x_3)=\{3\}$. It follows that $x_7\in A_1\cap A_2\cap A_3$ since
$x_7x_i\in E(D^*)$ for $i=1,2,3$. Thus, $A(x_7)=\{1,2,3\}$ and
$B(x_7)=\{4\}$. Note that $x_6x_i\in E(D^*)$ for $i=1,2,7$.  {Now}
$x_6\in A_1\cap A_2\cap A_4$ and therefore
$A(x_6)=A(x_3)=\{1,2,4\}$, contradicting Claim 1. {By} Claim 3, we have
$|A(x_3)|=|B(x_3)|=2$. Similarly, we obtain
$|A(x_7)|=|B(x_7)|=2$. Moreover, since $x_7x_i\in E(D^*)$ for
$i=1,2$, we have $A(x_7)=\{1,2\}$ and $B(x_7)=\{3,4\}$.

Note that $x_6x_7\in E(D^*)$.  {Hence $A(x_6)\cap B(x_7)\neq \emptyset$,
and} we assume without loss of generality that $x_6\in A_3$. This
implies that $A(x_6)=\{1,2,3\}$ and $B(x_6)=\{4\}$. Thus,
$B(x_4)=\{3\}$ or $B(x_5)=\{3\}$ since $|U|=4$. If $B(x_5)=\{3\}$,
then $A(x_3)=A(x_4)=\{3,4\}$, contradicting Claim 1. Therefore,
$B(x_4)=\{3\}$.

Now we have $B(x_1)=\{1\}$, $B(x_2)=\{2\}$, $B(x_4)=\{3\}$,
$B(x_6)=\{4\}$. It follows that $A(x_3)=\{3,4\}$, $A(x_5)=\{1,4\}$,
$A(x_7)=\{1,2\}$, and so $A(x_3)\cap A(x_5)\cap A(x_7)=\emptyset$.
Thus, the claim holds.
\end{proof}

\begin{claim}
$|A(z_{p,q,r})|=2$.
\end{claim}

\begin{proof}
It follows from $d^-_{D^*}(z_{p,q,r})\geq 1$ and
$d^+_{D^*}(z_{p,q,r})\geq 1$ that $1\leq |A(z_{p,q,r})|\leq 3$.
If $|A(z_{p,q,r})|=1$,  {then the edges from $z_{p,q,r}$
to $V(D_2)$ require $D_2$ to be covered by three  {cuts}},
contradicting Claim 2.
If $|A(z_{p,q,r})|=3$, then $|B(z_{p,q,r})|=1$. Since
$x_{p}z_{p,q,r},x_{q}z_{p,q,r},x_{r}z_{p,q,r}\in E(D^*)$, we have
$A(x_p)\cap A(x_q)\cap A(x_r)\neq \emptyset$, contradicting Claim 4.
Thus, $|A(z_{p,q,r})|=2$.
\end{proof}

Let $y_{p'},y_{q'},y_{r'}$ be the three vertices in $D_2$ such that
$|A(y_{p'})|=|A(y_{q'})|=|A(y_{r'})|=2$. Since there is a one-to-one
correspondence between the subsets of $\{1,2,3,4\}$ with cardinality
2 and the vertices in $\{x_p,x_q,x_r,y_{p'},y_{q'},y_{r'}\}$, we
have $A(z_{p,q,r})=A(v)$ for some $v\in
\{x_p,x_q,x_r,y_{p'},y_{q'},y_{r'}\}$, contradicting Claim 1.

Therefore, $D^*$ cannot be covered by four  {cuts}.

\section{Preliminaries}

Let $I=\{1,2,3,4,5\}$. We denote the ten subsets of $I$ with
cardinality 2 by  {$S_1,\ldots,S_{10}$} {(See Table 1)} and use them
to represent ten  {distinct} colors in the following.
 { We define a graph (with loops) on the colors by saying that
colors $S_i$ and $S_j$ are adjacent if and only if they share at least one element of $I$,
and say such two colors neighbor each other.}
Now we establish some properties of these ten colors
 {that} will be essential to our proof of Theorem 4.

\begin{center}

\begin{tabular}{c|c|c|c|c|c|c|c|c|c|c}
  \hline \hline
  $i$ & 1 & 2 & 3 & 4 & 5 & 6 & 7 & 8 & 9 & 10 \\
  \hline
  $S_i$ & $\{1,2\}$ & $\{1,3\}$ & $\{1,4\}$ & $\{1,5\}$ & $\{2,3\}$ & $\{2,4\}$ & $\{2,5\}$ & $\{3,4\}$ & $\{3,5\}$ & $\{4,5\}$ \\
  \hline \hline
\end{tabular}\\[3mm]

\footnotesize {Table 1: Subsets of $I$ with cardinality 2.}

\end{center}

\begin{prop}
 {
Colors $S_1,\ldots,S_{10}$ satisfy the following properties:\\
$(1)$ Each color neighbors itself and six other colors.\\
$(2)$ Any two colors have four common neighboring colors.\\
$(3)$ Any three colors have two common neighboring colors.\\
$(4)$ For any two pairs of distinct colors, some color in one pair neighbors some color in the other.
}
\end{prop}

\begin{proof}
{
$(1)$ Colors not neighboring $S_i$ consist of two of the three elements of $I$ not in $S_i$.}

{$(2)$ By (1), we may consider distinct colors $S_i$ and $S_j$.
If they are disjoint, then there are four ways to pick an element from each to form
a common neighboring color. Otherwise, the four colors containing their common element
are common neighboring colors (as is their symmetric difference).}

{$(3)$ Three colors cannot be pairwise disjoint.
If they have one common element,
then the four colors containing their common element are common neighboring colors.
Otherwise, the two colors containing the common element of some two adjacent colors and
one element in the remaining color are common neighboring colors.}

{$(4)$ If not, then there will be at least six different elements,
since two distinct colors contain at least three different elements. }
\end{proof}

For a digraph $D$, we color its vertices with the ten colors
 {$S_1,\ldots,S_{10}$} and use $c(v)$ to denote the color that has
been assigned to the vertex $v$ of $D$. For a
bipartition $(X,Y)$ of $V(D)$, let $X_i=\{x\in X\colon\,c(x)=S_i\}$ and
$Y_i=\{y\in Y\colon\,c(y)=S_i\}$. Our proof of Theorem 4 is heavily based
on the following lemma.

\begin{lemma}
{Let $D$ be a digraph, $(X,Y)$ a bipartition of $V(D)$ and $uv$ an edge of $D$}. If
there is a coloring of $V(D)$ with the ten colors
 {$S_1,\ldots,S_{10}$} {such that},\\
$(1)$ $c(u)$ and $c(v)$ are  {distinct} if $u,v\in X$ or $u,v\in Y$; and\\
$(2)$ $c(u)$ and $c(v)$ are  {adjacent if $u\in Y$ and $v\in X$,\\
then $D$ {can be} covered by the cuts $E(A_1,B_1),\ldots,E(A_5,B_5)$ defined by }
$$
A_i=(\bigcup_{i\notin S_j} X_j)\cup (\bigcup_{i\in S_j} Y_j),
B_i=(\bigcup_{i\in S_j} X_j)\cup (\bigcup_{i\notin S_j} Y_j).
$$
\end{lemma}

The following table describes the five  {cuts} of Lemma 1.

\begin{center}

\begin{tabular}{c|c|c}
  \hline \hline
  $i$ & $A_i$ & $B_i$ \\
  \hline
  1 & $X_5,X_6,X_7,X_8,X_9,X_{10},Y_1,Y_2,Y_3,Y_4$ & $X_1,X_2,X_3,X_4,Y_5,Y_6,Y_7,Y_8,Y_9,Y_{10}$ \\
  2 & $X_2,X_3,X_4,X_8,X_9,X_{10},Y_1,Y_5,Y_6,Y_7$ & $X_1,X_5,X_6,X_7,Y_2,Y_3,Y_4,Y_8,Y_9,Y_{10}$ \\
  3 & $X_1,X_3,X_4,X_6,X_7,X_{10},Y_2,Y_5,Y_8,Y_9$ & $X_2,X_5,X_8,X_9,Y_1,Y_3,Y_4,Y_6,Y_7,Y_{10}$ \\
  4 & $X_1,X_2,X_4,X_5,X_7,X_9,Y_3,Y_6,Y_8,Y_{10}$ & $X_3,X_6,X_8,X_{10},Y_1,Y_2,Y_4,Y_5,Y_7,Y_9$ \\
  5 & $X_1,X_2,X_3,X_5,X_6,X_8,Y_4,Y_7,Y_9,Y_{10}$ & $X_4,X_7,X_9,X_{10},Y_1,Y_2,Y_3,Y_5,Y_6,Y_8$ \\
  \hline \hline
\end{tabular}\\  [3mm]

\footnotesize {Table 2: The five  {cuts} $E(A_i,B_i)$,
$i=1,\ldots,5$.}

\end{center}

\begin{proof}
 {For $uv\in E(D)$ with $S_k=c(u)$ and $S_l=c(v)$,
we show that $uv\in E(A_i,B_i)$ for some $i$.}

 {First consider} $u\in X$ and $v\in Y$. Since $S_k\cup S_l$
contains at most four elements, there exists $i\in I$ such that
$i\notin S_k\cup S_l$.  {Now} $u\in A_i$ and $v\in B_i$ by the
 {definitions of $A_i$ and $B_i$.}

{If $u,v\in X$, then $S_k$ and $S_l$ are  {distinct}, and there exists
$i\in I$ such that $i\notin S_k$ and $i\in S_l$. If $u,v\in Y$, then $S_k$ and $S_l$ are  {distinct}, and there exists
$i\in I$ such that $i\in S_k$ and $i\notin S_l$. If $u\in Y$, $v\in X$, then $S_k$ and $S_l$ are  {adjacent}, and there
exists $i\in I$ such that $i\in S_k\cap S_l$. In all these cases, we have $u\in A_i$ and $v\in B_i$.}
\end{proof}

 {For} a digraph $D$ and a bipartition
$(X,Y)$ of $V(D)$, we say that a coloring of the vertices of $D$
with the  {ten colors $S_1,\ldots,S_{10}$ is a {\em good
coloring} for the bipartition $(X,Y)$ if} it satisfies the
conditions {of} Lemma 1.

\section{Proof of Theorem 4}

If every digraph in $\mathscr{D}(4,4)$ can be covered by  {five cuts,
then the same holds for} every digraph in $\mathscr{D}(3,3)$ since $\mathscr{D}(3,3)\subset \mathscr{D}(4,4)$,
 {hence} we only need to consider digraphs in $\mathscr{D}(4,4)$.

Let $D\in \mathscr{D}(4,4)$. We use $(X,Y)$ to denote an arbitrary bipartition
of $V(D)$ such that $d_D^-(x)\leq 4$ for every $x\in X$ and
$d_D^+(y)\leq 4$ for every $y\in Y$. By Lemma 1, it is sufficient to
show that there exists a good coloring of $D$ for the bipartition
$(X,Y)$.

If $|X|\leq 4$ and $|Y|\leq 4$, then we can use $S_1,S_2,S_3,S_4$
(all contain the element 1) to color the vertices of $D$  {so} that
any two vertices in $X$ ($Y$) receive  {distinct} colors. It follows
that  {$c(x)$ neighbors $c(y)$} for any $x\in X$ and $y\in Y$,
and so we get a good coloring of $D$ for the bipartition $(X,Y)$.

Suppose now that $\max \{|X|,|Y|\}\geq 5$ and  there exist no good
colorings of $D$ for the bipartition $(X,Y)$. We assume that $D$ is
chosen such that $|V(D)|$ is as small as possible.

\setcounter {claim}{0}
\begin{claim}
$d_{X}^-(x)=d_{X}^+(x)$ for every $x\in X$.
\end{claim}

\begin{proof}
 {If not, then} we have $d_{X}^-(x)>d_{X}^+(x)$ for some
$x\in X$. This implies that $d_{Y}^-(x)\leq 3$.
 {Now $D-\{x\}$} has a good coloring for the bipartition $(X\backslash
\{x\},Y)$ by the choice of $D$. Let $c_1,\ldots,c_t$ be the colors
of the vertices in $X$  {that neighbor} $x$.

Assume first that $d_{Y}^-(x)=0$. Since $d_{X}^+(x)<d_{X}^-(x)\leq
4$, we have $t\leq d_{X}(x)\leq 7$. It follows that there exists a
color $c\notin \{c_1,\ldots,c_t\}$.  {Assign $c$ to $x$ to obtain a
good coloring of $D$ for the bipartition $(X,Y)$. }

{Suppose that $d_{Y}^-(x)=1$, and choose $y\in Y$ such that $yx\in E(D)$.} Since
$d_{X}^+(x)<d_{X}^-(x)\leq 3$, we have $t\leq d_{X}(x)\leq 5$.
 {By Proposition 1 (1)}, there exists a color $c\notin \{c_i\colon\,i=1,\ldots,t\}$
 {that neighbors} $c(y)$.  {Assign $c$ to $x$ to obtain a good
coloring of $D$ for the bipartition $(X,Y)$. }

{Suppose that $d_{Y}^-(x)=2$, and choose $y_1,y_2\in Y$ such that $y_1x,y_2x\in
E(D)$.} Since $d_{X}^+(x)<d_{X}^-(x)\leq 2$, we have $t\leq
d_{X}(x)\leq 3$.  {By Proposition 1 (2)}, there exists a color $c\notin
\{c_i\colon\,i=1,\ldots,t\}$  {that neighbors} both $c(y_1)$ and $c(y_2)$.
 {Assign $c$ to $x$ to obtain a good coloring of $D$ for
the bipartition $(X,Y)$.}

{Suppose that $d_{Y}^-(x)=3$, and choose $y_1,y_2,y_3\in Y$ such that
$y_1x,y_2x,y_3x\in E(D)$.} Since $d_{X}^+(x)<d_{X}^-(x)\leq 1$, we
have $t\leq d_{X}(x)\leq 1$.  {By Proposition 1 (3)}, there exists a color
$c\notin \{c_i\colon\,i=1,\ldots,t\}$  {that neighbors} any one in
$\{c(y_1),c(y_2),c(y_3)\}$.  {Assign $c$ to $x$ to obtain a good
coloring of $D$ for the bipartition $(X,Y)$.}
\end{proof}

\begin{claim}
$d_{Y}^-(x)\geq 2$ for every $x\in X$.
\end{claim}

\begin{proof}
Assume that $d_{Y}^-(x)\leq 1$ for some $x\in X$.
 {Now $D-\{x\}$} has a good coloring for the bipartition $(X\backslash
\{x\},Y)$. Let $c_1,\ldots,c_t$ be the colors of the vertices in $X$
 {that neighbor} $x$.

Assume first that $d_{Y}^-(x)=0$. Since $d_{X}^+(x)=d_{X}^-(x)\leq 4$,
we have $t\leq d_{X}(x)\leq 8$. It follows that there exists a
color $c\notin \{c_1,\ldots,c_t\}$.  {Assign $c$ to $x$ to obtain a
good coloring of $D$ for the bipartition $(X,Y)$.}

If $d_{Y}^-(x)=1$, then choose $y\in Y$ such that $yx\in E(D)$. Since
$d_{X}^+(x)=d_{X}^-(x)\leq 3$, we have $t\leq d_{X}(x)\leq 6$.
 {By Proposition 1 (1)}, there exists a color $c\notin \{c_1,\ldots,c_t\}$
 {that neighbors} $c(y)$.  {Assign $c$ to $x$ to obtain a good
coloring of $D$ for the bipartition $(X,Y)$.}
\end{proof}

Similarly we can get

\begin{claim}
$d_{Y}^-(y)=d_{Y}^+(y)$ for every $y\in Y$.
\end{claim}

\begin{claim}
$d_{X}^+(y)\geq 2$ for every $y\in Y$.
\end{claim}

It follows from Claim 2 that there exists an edge $yx\in E(D)$,
where $x\in X$ and $y\in Y$. By the choice of $D$, there is a good
coloring of $D-\{x,y\}$ for the bipartition $(X\backslash
\{x\},Y\backslash \{y\})$. Let $c_1,\ldots,c_t$ be the colors of the
vertices in $X$  {that neighbor} $x$.

\begin{claim}
 {There exist two distinct colors $b_1,b_2$ that we can assign to
$x$ so that} the resulting coloring of $D-\{y\}$ is a good coloring
for the bipartition $(X,Y\backslash \{y\})$.
\end{claim}

\begin{proof}
Let $Y'=Y\backslash \{y\}$. Since $2\leq d^-_Y(x)\leq 4$ and $yx\in
E(D)$, we have $1\leq d^-_{Y'}(x)\leq 3$.

{If $d^-_{Y'}(x)=1$, then choose $y_1\in Y'$ such that $y_1x\in E(D)$.}
Since $d_{X}^+(x)=d_{X}^-(x)\leq 2$, we have $t\leq d_{X}(x)\leq 4$.
 {By Proposition 1 (1)}, there exist two  {distinct} colors $b_1,b_2\notin
\{c_i\colon\,i=1,\ldots,t\}$  {that both neighbor} $c(y_1)$. We
can get a good coloring of $D-\{y\}$ for the bipartition
$(X,Y\backslash \{y\})$ by assigning either $b_1$ or $b_2$ to $x$.

{If $d_{Y'}^-(x)=2$, then choose $y_1,y_2\in Y'$ such that $y_1x,y_2x\in
E(D)$.} Since $d_{X}^+(x)=d_{X}^-(x)\leq 1$, we have $t\leq d_{X}(x)\leq 2$.
 {By Proposition 1 (2)}, there exist two  {distinct} colors
$b_1,b_2\notin \{c_i\colon\,i=1,\ldots,t\}$  {that both neighbor}
any one in $\{c(y_1),c(y_2)\}$. We can get a good coloring of
$D-\{y\}$ for the bipartition $(X,Y\backslash \{y\})$ by assigning
either $b_1$ or $b_2$ to $x$.

{If $d_{Y'}^-(x)=3$, then choose $y_1,y_2,y_3\in Y'$ such that
$y_1x,y_2x,y_3x\in E(D)$.} It is clear that $d^+_X(x)=d^-_X(x)=0$.
 {By Proposition 1 (3)}, there exist two  {distinct} colors $b_1,b_2$
 {that both neighbor} any one in $\{c(y_1),c(y_2),c(y_3)\}$.
We can get a good coloring of $D-\{y\}$ for the bipartition
$(X,Y\backslash \{y\})$ by assigning either $b_1$ or $b_2$ to $x$.

Thus, the claim holds.
\end{proof}

Similarly,  {there exist two distinct colors $b_3,b_4$ that we can
assign to $y$ so that} the resulting coloring of $D-\{x\}$ is a
good coloring for the bipartition $(X\backslash \{x\},Y)$.
 {For the two pairs of distinct colors $b_1,b_2$
and $b_3,b_4$, by Proposition 1 (4), $b_i$ neighbors $b_j$ for
some $b_i\in \{b_1,b_2\}$ and some $b_j\in \{b_3,b_4\}$.}
Assign $b_i$ and $b_j$ to $x$ and $y$, respectively.
 {This produces} a good coloring
of $D$ for the bipartition $(X,Y)$, a contradiction.

The proof of Theorem 4 is complete.


\end{document}